\documentclass[12pt]{article}
\usepackage{latexsym}
\usepackage{amssymb}
\begin{document}
\newtheorem{theorem}{Theorem}[section]
\newtheorem{lemma}[theorem]{Lemma}
\newtheorem{proposition}[theorem]{Proposition}
\newtheorem{corollary}[theorem]{Corollary}
\newtheorem{definition}{Definition}
\newtheorem{question}{Question}
\newtheorem{conjecture}{Conjecture}
\newcommand{\F}{\ensuremath{\mathbb F}}
\newcommand{\N}{\mathcal N}
\newcommand{\R}{\mathcal R}
\newcommand{\Z}{\mathbb Z}
	
\title{Conjugacy classes of maximal cyclic subgroups of metacyclic $p$-groups}
\author{M. Bianchi, R.D. Camina, \&  Mark L. Lewis}
\date \
\maketitle
	
\begin{abstract} 
In this paper, we set $\eta (G)$ to be the number of conjugacy classes of maximal cyclic subgroups of a finite group $G$.  We compute $\eta (G)$ for all metacyclic $p$-groups.  We show that if $G$ is a metacyclic $p$-group of order $p^n$ that is not dihedral, generalized quaternion, or semi-dihedral, then $\eta (G) \ge n-2$, and we determine when equality holds.  \\[1ex]
\end{abstract}
	
	\maketitle
	
\section{Introduction}
	
Unless otherwise stated, all groups in this paper are finite, and we will follow standard notation from \cite{isstext}.  As in \cite{pre1} and \cite{pre2}, we set $\eta (G)$ to be the number of conjugacy classes of maximal cyclic subgroups of a group $G$.  For $p = 2$, we have  that $\eta (G) = 3$ when $G$ is a dihedral $2$-group, a generalized quaternion $2$-group, or a semi-dihedral group.  In \cite{pre}, the second and third authors along with Yiftach Barnea and Mikhail Ershov have shown that for every prime $p \ge 5$ there are infinitely many $p$-groups with $\eta = p + 2$ and for $p =3$ there are infinitely many $3$-groups with $\eta = 9$. This answers negatively Question 5.0.9 from \cite{von} which asked whether $\eta(G)$ grows with the order of $G$ when $G$ is a $p$-group and $p$ is odd.

On the other hand, it is rare for this to occur.  	Indeed, the only $2$-groups (in fact the only $p$-groups) that have $\eta = 3$ are the Klein $4$-group, the dihedral groups, the generalized quaternion groups, and the semi-dihedral groups.  To see this, we know that $\eta (G) \ge \eta (G/G')$ (see \cite{pre1}), and for $p$-groups $\eta (G/G') \ge p+1$ when $G/G'$ is not cyclic (see \cite{pre2}).  Thus, $\eta = 3$ can only occur when $p = 2$.  Also, in \cite{pre2}, we show that $\eta (G/G') = 3$ if and only if $G/G' \cong C_2 \times C_2$.  It is well known that if $G$ is a $2$-group of order at least $8$ and $|G:G'| = 4$, then $G$ is either dihedral, generalized quaternion, or semi-dihedral.  (See Problem 6B.8 of \cite{isstext}.)

Now, dihedral groups, generalized quaternion groups, and semi-dihedral groups are examples of metacyclic groups.  I.e., groups $G$ with a normal subgroup $N$ so that $N$ and $G/N$ are both cyclic groups.  
This motivated us to investigate the invariant $\eta$ for all metacyclic $p$-groups. Indeed this project began before the results of \cite{pre} were known and
we were originally curious as to whether  we would find another family of metacyclic $p$-groups with fixed $\eta$. However, we prove the following:
	
\begin{theorem} \label{main1}
Let $G$ be a metacyclic $p$-group of order $p^n$ that is not a dihedral group, generalized quaternion group, or semi-dihedral group.  Then $\eta (G) \ge n - 2$.
\end{theorem}
	
In fact, we compute $\eta (G)$ for every metacyclic $p$-group $G$.  Thus, we list the metacyclic $p$-groups where equality occurs in Theorem \ref{main1}.  King in \cite{king} gave a description of all metacyclic $p$-groups.  We will give this description of these groups in Section \ref{secn:meta}.  In particular, King divided the metacyclic $p$-groups into two families of groups which he called {\it positive type} and {\it negative type}.  The negative type groups only occur when $p = 2$, so if $p$ is an odd prime, then all of the metacyclic $p$-groups are of positive type.  We have the following result for the metacyclic groups of positive type.
	
\begin{theorem} \label{main2}
Let $G$ be a metacyclic group of positive type.  Then $\eta (G) = \eta (G/G')$.
\end{theorem}
	
We note that Rog\'erio in \cite{rogerio} has a formula to compute $\eta (A)$ for an abelian group $A$.  His formula involves the Euler $\phi$-function and a second number theoretic function.  When $G$ is a metacyclic abelian $p$-group, we prove in \cite{pre2} a formula for $\eta (G)$ that is only in terms of the sizes of the direct factors of $G$.  Notice in Theorem \ref{main2} that $G/G'$ will be a metacyclic abelian $p$-group, and so, our formula will compute $\eta (G/G')$ and hence, $\eta (G)$.
	
When $G$ is a metacyclic $p$-group of negative type, it is not usually the case that $\eta (G)$ and $\eta (G/G')$ are equal.  However, we will find that there usually is a proper quotient whose value of $\eta$ equals $\eta (G)$.  We will also see for most metacyclic groups of negative type that the formula for $\eta$ is dependent on the formula for $\eta$ that we found for the metacyclic abelian $p$-groups.

The authors would like to thank Emanuele Pacifici for a number of helpful conversations while working on this paper.	
	
\section{Preliminaries}
	
In our preprint \cite{pre1}, we prove two results that we need in this paper.  The first is a criteria for determining when the quotient of a $p$-group $G$ has the same value for $\eta$ as $\eta (G)$.  Given a prime $p$, we set $G^{\{p\}} = \{ g^p \mid g \in G \}$.  I.e., $G^{\{p\}}$ is the set of $p$-th powers in $G$.

\begin{theorem}\label{quot} 
Let $N$ be a normal subgroup of the $p$-group $G$. Then $\eta(G/N)  \leq  \eta(G).$  Furthermore, $\eta(G/N) = \eta(G)$ if and only if $N\subseteq G^{\{p\}}$ and for all $x \in G \setminus G^{\{p\}}$ every element of $xN$ is conjugate to a generator of $\langle x \rangle$.  In particular, if $\eta(G/N) = \eta(G)$, then $G^{\{p\}}$ is a union of $N$-cosets and $G^{\{p\}}N = G^{\{p\}}$. 
\end{theorem}

This second Proposition relates $\eta (G)$ to the number of $G$-orbits of maximal cyclic subgroups of a normal subgroup.

\begin{proposition}\label{centre}  
Let $N$ be a normal subgroup of a group $G$ and let $\eta^*(N)$ be the number of $G$-orbits on the $N$-conjugacy classes of maximal cyclic subgroups of $N$.  Then $\eta (G) \geq \eta^*(N)$.  In particular,
\begin{enumerate}	
\item[(i)] if $N$ is central in $G$, then $\eta(G) \geq \eta(N)$.
\item[(ii)] if $|G:N| = k$, then $\eta (G) \geq \eta(N)/k$.
\end{enumerate}
\end{proposition}

Let $p$ be a prime, and let $a$ and $b$ be positive integers.  We take $k = {\rm max} (a,b)$ and $l = {\rm min} (a,b)$.  We set $g_p (a,b) = p^{(l-1)}((k-l) (p-1) + p +  1)$.  In \cite{pre2}, we prove the following lemma.

\begin{lemma}\label{two}
If $p$ is a prime and $a$ and $b$ are positive integers so that $G = C_{p^a} \times C_{p^b}$, then $g_p (a,b) = \eta (G)$.
\end{lemma}

We close this section with an easy lemma that computes $g_2$ for small values and gives a lower bound for larger values.  We remark that when $p = 2$, this function is much easier to work with.

\begin{lemma}\label{g2 comp}
Suppose $k \ge l$.  Then the following hold:
\begin{enumerate}
\item If $l = 1$, then $g_2 (k,1) = k+2$.
\item If $l = 2$, then $g_2 (k,2) = 2(k+1)$.		
\item If $l = 3$, then $g_2 (k,3) = 4k$.
\item If $l \ge 4$, then $g_2 (k,l) \ge 4k + 2l$. 
\end{enumerate}
\end{lemma}

\noindent
{\bf Proof.}
We have $g_2 (a,b) = g_2 (k,l) = 2^{l-1}(k-l+3)$.  Conclusions (1), (2), and (3) are immediate.  We focus on (4).  Begin with $g_2 (4,4) =24$; so the result holds for $g_2 (4,4)$.  Next, $g_2 (l,l) - 6l = 3 \cdot 2^{l-1} - 6l$ is clearly increasing when $l \ge 3$.  Thus, we have $g_2 (l,l) \ge 4l + 2l$ when $l \ge 3$.  Let $k = l + m$ for $m \ge 0$.  Then $g_2 (k,l) = g_2 (l+m,m) = 2^{l-1}(m+3)$ and $4k + 2l = 4(l+m) + 2l = 6l + 4m$.  Fixing $l \ge 4$, we note that $2^{l-1}(m+3) - 6l - 4m$ will be an increasing function in $m$.  We conclude that $g_2 (k,l) \ge 4k + 2l$ for $l \ge 4$.
$\Box$\\

\section{Metacyclic $p$-Groups} \label{secn:meta}

For the rest of the paper, we will focus on metacyclic $p$-groups.  A finite metacyclic $p$-group can be described as follows.  This description is taken from \cite{king}, 
$$G_p(\alpha, \beta, \epsilon, \delta, \pm ) = \langle x,y \mid x^{p^{\alpha}}=1, y^{p^{\beta}}= x^{p^{\alpha - \epsilon}}, x^y = x^r \rangle$$
where $r = p^{\alpha - \delta} +1$ (positive type) or $r = p^{\alpha - \delta} -1$ (negative type).
The integers $\alpha$, $\beta$, $\delta$, $\epsilon$ satisfy $\alpha, \beta >0$ and $\delta, \epsilon$
nonnegative, furthermore $\delta \leq {\rm min} \{\alpha -1, \beta \}$ and $\delta + \epsilon \leq \alpha$.  When $G$ has negative type, only $\epsilon = 0$ or $1$ occur.  For $p$ odd
$$G \cong G_p(\alpha, \beta, \epsilon, \delta, +).$$  
In other words, the negative type only occurs when $p = 2$; when $p$ is odd, only the positive type occurs.  Metacyclic $2$-groups can be of either positive type or negative type.  We note that dihedral, semi-dihedral and generalized quaternion groups are all of negative type. 

If $p=2$, then in addition $\alpha - \delta >1$ and
$$G \cong G_2(\alpha, \beta, \epsilon, \delta, +)\;\;{\rm or}\;\;G \cong G_2(\alpha, \beta, \epsilon, \delta, -).$$
Note, the above presentation does not guarantee nonisomorphic groups for different parameters (see \cite{beuerle}). However, the parameters do determine some structural information about $G$. For example, $|G| = p^{\alpha + \beta}$ and $G' = \langle x^{p^{\alpha - \delta}} \rangle$ if $G$ is of positive type and $G' = \langle x^2 \rangle$ if $G$ is of negative type. All elements of $G$ can be written as $y^bx^a$ for some integers $a$ and $b$.  Also if $G$ is of positive type then  $Z(G) = \langle x^{p^{\delta}}, y^{p^{\delta}} \rangle$ and $|Z(G)| = p^{\alpha + \beta - 2 \delta}$, if $G$ is of negative type $Z(G) = \langle x^{2^{\alpha -1}}, y^{2^{\max \{1, \delta\}}} \rangle$, \cite[Prop. 2.5]{beuerle}.  Note that if $G$ is of positive type and $\delta = 0$, then $G$ will be abelian.  

As we mentioned above, the dihedral groups, the generalized quaternion groups, and the semi-dihedral groups are the only $p$-groups $G$ that satisfy $\eta(G) = 3$.  These are also precisely the $2$-groups of maximal class.  We have also mentioned that they are metacyclic.  In terms of our notation, the dihedral groups are $G_2 (\alpha,1,0,0,-)$, the generalized quaternion groups are $G_2 (\alpha,1,1,0,-)$, and the semi-dihedral groups are $G_2 (\alpha,1,0,1,-)$.

For Lemmas \ref{quo1} and \ref{conjugacy}, we are writing $G_p (\alpha,\beta,\epsilon,\delta, \pm )$ as $G_p (\alpha, \beta,\epsilon,\delta, \gamma)$ where we take $\gamma = +$ when $G$ is of positive type and $\gamma = -$ when $G$ is of negative type.   We consider quotients of $G$.  Note that this lemma would not be well defined if $\delta = 0$ and would not say anything if $\delta = 1$.

\begin{lemma} \label{quo1}
Suppose $G$ is $G_p (\alpha, \beta, \epsilon, \delta, \gamma)$ with $\delta \ge 2$.  Then  $N = \langle x^{p^{\alpha - \delta + 1}} \rangle$ is a normal subgroup of $G$ and $G/N$ is isomorphic to 
$$G_p (\alpha-\delta+1, \beta, (\epsilon - \delta +1)^*, 1, \gamma )$$ where $(\epsilon - \delta +1)^* = \epsilon - \delta +1$ when $\epsilon \ge \delta - 1$ and $(\epsilon - \delta +1)^* = 0$ when $\epsilon < \delta -1$.
\end{lemma}

\noindent
{\bf Proof.}
Set $Z = \langle x^{p^{\alpha -1}} \rangle \le Z (G)$.  We first prove that $G/Z$ is isomorphic to $G_p (\alpha -1, \beta, \epsilon - 1, \delta - 1, \gamma)$ when $\epsilon \ge 1$ and $G_p (\alpha -1, \beta, 0, \delta - 1, \gamma)$ when $\epsilon = 0$. We know that $G/Z = \langle xZ, yZ \rangle$ where $xZ$ has order $p^{\alpha - 1}$.  Observe that $(yZ)^{p^\beta} = y^{p^\beta} Z = x^{p^{\alpha -\epsilon}}Z$.  When $\epsilon \ge 1$, we have $$x^{p^{\alpha -\epsilon}}Z = x^{p^{(\alpha - 1) - (\epsilon -1)}}Z$$ and when $\epsilon = 0$, we have $$x^{p^{\alpha -\epsilon}}Z = x^{p^\alpha}Z = Z.$$  Also, $$(xZ)^{yZ} = x^y Z = x^{p^{\alpha - \delta} + \gamma} Z = x^{p^{(\alpha - 1) - (\delta - 1)}+ \gamma}Z.$$  Hence, $G/Z$ satisfies the hypotheses for  $G_p (\alpha -1, \beta, \epsilon - 1, \delta - 1, \gamma)$ when $\epsilon \ge 1$ and $G_p (\alpha -1, \beta, 0, \delta - 1, \gamma)$ when $\epsilon = 0$.

We know that $X = \langle x \rangle$ is a cyclic, normal subgroup of $G$.  Observe that $N$ is contained in $X$ and so is characteristic.  This implies that $N$ is normal in $G$.   Observe that $Z \le N$ and we have shown that
$G/Z \cong G_p (\alpha -1, \beta, \epsilon - 1, \delta - 1, \gamma)$ or $G_p (\alpha -1, \beta, 0, \delta - 1, \gamma)$.  If $\delta = 2$, then $N = Z$, and we have the desired result.  Otherwise, we have $\delta \ge 3$.  Using induction, we have $G/N \cong (G/Z)/(N/Z)$ is isomorphic to either 
$$
G_p((\alpha-1)-(\delta - 1)+1, \beta, (\epsilon - 1) - (\delta - 1) +1, 1, \gamma ) 
\cong  G_p (\ \alpha-\delta+1, \beta, \epsilon - \delta + 1, 1, \gamma )
$$
or 
$$ G_p((\alpha-1)-(\delta - 1)+1, \beta, 0, 1, \gamma ) \cong G_p ( \alpha-\delta+1, \beta, 0, 1, \gamma ). ~~\Box
$$
\\

We consider the metacyclic groups of positive type and use Theorem \ref{quot} and Lemma \ref{two}. Thus, we first analyze $G/G'$. 

\begin{lemma}\label{ten} 
Suppose $G = G_p(\alpha, \beta, \epsilon, \delta, +)$.
\begin{enumerate}
\item [(i)]  If $\delta \geq \epsilon$ or $\delta < \epsilon$ and $\alpha \geq \beta + \epsilon$, then $G/G' =  C_{p^{\alpha - \delta}} \times C_{p^{\beta}}$.
\item [(ii)] If $\delta < \epsilon$ and $\alpha < \beta + \epsilon$, then $G/G' = C_{p^{\alpha - \epsilon}}  \times C_{p^{\beta + \epsilon -\delta}}$.
\end{enumerate}
\end{lemma}

\noindent
{\bf Proof.}
Now $G' =  \langle x^{p^{\alpha - \delta}} \rangle$, so $|G'| = p^{\delta}$. Also  $|G| = p^{\alpha + \beta}$, 
so $|G:G'| = p^{\alpha +\beta - \delta}.$  

If $\delta \geq \epsilon$, then $\langle y \rangle \cap G' =  \langle x^{p^{\alpha-\epsilon}} \rangle =  \langle x \rangle \cap \langle y \rangle$.  We see that $xG'$ has order $p^{\alpha - \delta}$, and $yG'$ has order $p^\beta$ 
and $G/G' =  \langle xG' \rangle  \times  \langle yG' \rangle$ yielding the desired result.

Now suppose $\delta < \epsilon$.  In this case, we see that $G' <  \langle x^{p^{\alpha - \epsilon}} \rangle  =  \langle x \rangle  \cap \langle y \rangle$. We see that $xG'$ has order $p^{\alpha - \delta}$ and $yG'$ has order $p^{\beta + 	\epsilon - \delta}$.  Since $G' <  \langle x \rangle  \cap  \langle y \rangle$, we do not have that $G/G'$ is a direct product of $\langle xG' \rangle$ and $ \langle yG' \rangle$.  We see that $G/G'$ is abelian and generated by $xG'$ and $yG'$, so every element of $G/G'$ has order $\le {\rm max} \{p^{\alpha - \delta}, p^{\beta+\epsilon-\delta} \}$.  If $\alpha \geq \beta +\epsilon$, then $\alpha - \delta \geq \beta + \epsilon - \delta$.  In this case, $xG'$ has the largest order of any element in $G/G'$, and so we get $G/G' =  C_{p^{\alpha - \delta}} \times C_{p^\beta}$ since $|G/G'| = p^{\alpha + \beta - \delta}$.  On the other hand, if $\alpha < \beta + \epsilon$, then $\alpha - \delta < \beta + \epsilon - \delta$.  In this case, $yG'$ has the largest order of any element in $G/G'$ and we get $G/G'  = C_{p^{\alpha - \epsilon}} \times C_{p^{\beta + \epsilon - \delta}}$.
$\Box$\\

Given an element $g \in G$, we write ${\rm cl} (g)$ to denote the conjugacy class of $g$ in $G$.

\begin{lemma} \label{conjugacy}
	Let $G = G_p (\alpha, \beta, \epsilon, \delta, \gamma )$.  If $g = y^{pl+a}x^m$ for integers $l$, $m$, and $a$ so that $a \in \{ 1, \dots, p-1\}$, then ${\rm cl} (g) = g G'$.
\end{lemma}

\noindent
{\bf Proof.}
We first claim that $G = \langle x, g \rangle$.  We know that $G = \langle x, y \rangle$.  Obviously, $\langle x, g \rangle \le G$.  Observe that $y^{pl + a} = gx^{-m} \in \langle x, g \rangle$.  Since the order of $y$ is a power of $p$, this implies that $y \in \langle x, g \rangle$.  We conclude that $G = \langle x, y \rangle \le \langle x, g \rangle \le G$.  This proves the claim.  

Because $\langle x \rangle$ is normal in $G$, we obtain $G = \langle x \rangle \langle g \rangle$.  Observe that $\langle g \rangle \le C_G (g)$. By Dedekind's lemma (see Lemma X.3 on page 328 of \cite{isstext}), it follows that $C_G (g) = (C_G(g) \cap \langle x \rangle ) \langle g \rangle = C_{\langle x \rangle} (g) \langle g \rangle $.  Since $x$ centralizes $x^m$, we have 
$$C_{\langle x \rangle} (g) = C_{\langle x \rangle} (y^{pl+a}x^m) = C_{\langle x \rangle} (y^{pl+a}) = C_{\langle x \rangle} (y) = \langle x^{p^t} \rangle,$$ where $t = \delta$ if $\gamma = +$ and $t = \alpha - 1$ when $\gamma = -$.
We see that $C_G (g) = \langle g, x^{p^t} \rangle$.  We deduce that 
$$|G:C_G (g)| = |\langle x \rangle:\langle x^{p^t} \rangle| = p^t = |G'|.$$    Since ${\rm cl} (g) \subseteq gG'$, we conclude that ${\rm cl} (g) = g G'$.  
$\Box$ \\

Given a group $G$ and a prime $p$, we define $G^p = \langle G^{\{p\}} \rangle$.  I.e., $G^p$ is the subgroup generated by $G^{\{p\}}$.  In a similar fashion, we define $G^4 = \langle g^4 \mid g \in G \rangle$.  Following the literature, we say that a finite $p$-group $G$ is {\it powerful} if (i) $G' \le G^p$ when $p$ is odd and (ii) $G' \le G^4$ when $p = 2$. If $G$ is a powerful $p$-group, then it is known that $G^p = G^{\{p\}}$, i.e. the set of $p$-powers of elements of $G$ is equal to the subgroup the $p$-powers generate.  (See Section 2 of \cite{powerful} and in particular Propostion 2.6 of that citation.)

We claim that metacyclic $p$-groups of positive type are powerful. Let $G$ be $G_p (\alpha, \beta, \epsilon, \delta, +)$, then $G' = \langle x^{p^{\alpha - \delta}} \rangle.$  As $\alpha - \delta \ge 1$ it follows immediately that $G$ is powerful when $p$ is odd. For $p = 2$, we note that $\alpha - \delta \ge 2$ so again we have that $G$ is powerful.  

When $G$ is of positive type, we extend Lemma \ref{conjugacy}.

\begin{lemma} \label{unnamed}
Let $G = G_p (\alpha, \beta, \epsilon, \delta, +)$ and $g \in G \setminus G^{\{p\}}$. Then ${\rm cl}(g) =
gG'$.
\end{lemma}

\noindent
{\bf Proof.} 
Let $g \in G$ then $g = y^nx^m$ for some integers $n$ and $m$. As $G$ is
powerful, it follows that if $g \in G \setminus G^{\{p\}}$, then $g \not\in G^p$, and thus, one of $n$ and $m$
is not divisible by $p$.  When $n$ is not divisible by $p$, we obtain the conclusion by Lemma \ref{conjugacy}.

We now suppose that $g = y^nx^m$ where $m$ is not divisible by $p$.  We want to prove that ${\rm cl} (g) = gG'$.  We know that ${\rm cl} (g) \subseteq gG'$.  It suffices to prove that $|{\rm cl} (g)| \ge |gG'| = |G'| = p^\delta$.  On the other hand, we know that $y$ acts as an automorphism of order $p^\delta$ on $\langle x \rangle$, so $x$ has $p^\delta$ distinct images under powers of $y$.  Thus, if $1 \le a,b \le p^\delta$, then $x^{y^a} = x^{y^b}$ if and only if $a = b$.  Since $m$ is coprime to $p$, we see that $(x^{y^a})^m = (x^{y^b})^m$ if and only if $a = b$.  Hence, we have that $g^{y^a} = g^{y^b}$ if and only if $(y^nx^m)^{y^a} = (y^nx^m)^{y^b}$ and this occurs if and only if $a = b$.  We deduce that $g$ has at least $p^\delta$ distinct conjugates under $\langle y \rangle$ and so $|{\rm cl} (g)| \ge p^{\delta}$ as desired.  This proves the lemma.
$\Box$\\

We now prove that if $G$ is metacyclic of positive type, then $\eta (G) = \eta (G/G')$.  Combining this fact with Lemmas \ref{two} and \ref{ten}, we are able to compute $\eta (G)$ for all primes $p$.

\begin{corollary} \label{positive quo}
Suppose $G$ is $G_p(\alpha, \beta, \epsilon, \delta, +)$. Then $\eta(G) = \eta(G/G')$.	
\end{corollary}

\noindent
{\bf Proof.} 
As $G$ is powerful, by Theorem \ref{quot}, we need to show that for all
$g \in G \setminus G^{\{p\}}$ every element of $gG'$ is conjugate to a generator of 
$\langle g \rangle$, this follows from Lemma \ref{unnamed}.
$\Box$\\

For the record, we explicitly record the value of $\eta (G)$ when $G$ is a metacyclic group of positive type.

\begin{corollary} \label{explicit}
Suppose $G$ is $G_p(\alpha, \beta, \epsilon, \delta, +)$.
\begin{enumerate}
\item[(i)] If $\delta \ge \epsilon$ or $\delta < \epsilon$ and $\alpha \ge \beta + \epsilon$, then $\eta (G) = g_p (\alpha -\delta, \beta)$.
	\begin{enumerate}
	\item If $\beta \le \alpha - \delta$, then $\eta (G) = p^{\beta - 1} ((\alpha - \delta - \beta)(p-1) + p +1)$.
	\item If $\beta > \alpha - \delta$, then $\eta (G) = p^{\alpha - \delta - 1} ((\beta - \alpha + \delta)(p-1) + p + 1)$.
	\end{enumerate}
\item[(ii)] If $\delta < \epsilon$ and $\alpha < \beta + \epsilon$, then $\eta (G) = g_p (\alpha - \epsilon, \beta+\epsilon - \delta)= p^{\alpha -\epsilon - 1} ((\beta - \alpha + 2\epsilon - \delta)(p-1) + p + 1)$.

\end{enumerate}	
\end{corollary}

\noindent
{\bf Proof.}
Using Corollary \ref{positive quo}, we have $\eta (G) = \eta (G/G')$.  If $\delta \ge \epsilon$ or $\delta < \epsilon$ and $\alpha \ge \beta + \epsilon$, then in view of Lemma \ref{ten}, we see that $G/G' =  C_{p^{\alpha - \delta}} \times C_{p^{\beta}}$ and $\eta (G) = g_p (\alpha - \delta, \beta)$.  The remainder of (i) follows from the definition of $g_p$.  Suppose $\delta < \epsilon$ and $\alpha < \beta + \epsilon$.  Applying Lemma \ref{ten}, we see that $G/G' = C_{p^{\alpha - \epsilon}}  \times C_{p^{\beta + \epsilon -\delta}}$.  Observe that $\alpha < \beta + \epsilon$ yields $\alpha - \epsilon < \beta < \beta + \epsilon - \delta$ as we are assuming $\delta < \epsilon$.  In light of the definition of $g_p$, we obtain conclusion (ii).
$\Box$\\

When $G$ is metacyclic of positive type, we show that $\eta (G) \ge \alpha + \beta$.

\begin{corollary}
If $G$ is $G_p (\alpha,\beta,\epsilon,\delta,+)$, then $\eta (G) \ge \alpha + \beta$.
\end{corollary}

\noindent
{\bf Proof.}
We consider separately the cases given in Corollary \ref{explicit}. We use the fact that $2^{\beta - 1} \geq \beta$ for $\beta$ a positive integer. First, (i)(a),
where $\alpha - \delta \geq \beta$,
\begin{eqnarray*}
	\eta(G) & = & p^{\beta -1}((\alpha - \delta - \beta)(p-1) + p + 1) \\
	& \geq & 2^{\beta -1} (\alpha - \delta - \beta + 3 )\\
	& \geq & \alpha - \delta - \beta + 3\beta\\
	& \geq & \alpha + \beta
\end{eqnarray*}
since $\beta \geq \delta$.

Now, case (i)(b), so $\alpha - \delta < \beta$. First assume $\alpha - \delta > 1$, then 
\begin{eqnarray*}
	\eta(G) & = & p^{\alpha - \delta - 1}((\beta - \alpha + \delta)(p-1) + p+1)\\
	& \geq & 2^{\alpha - \delta -1}(\beta - \alpha + \delta + 3)\\
	& \geq & 2(\beta - \alpha + \delta) + 3(\alpha - \delta)\\
	& = & 2\beta +  (\alpha - \delta)\\
	& \geq & \beta + \alpha + (\beta - \delta) \\
	& \geq & \beta + \alpha
\end{eqnarray*}
since $\beta \geq \delta$. If $\alpha - \delta =1$ then $p \geq 3$, also note $\alpha = 1 + \delta \leq 1 + \beta$. So, we have 
$$\eta(G) \geq 2(\beta - \alpha + \delta) + 4  = 2\beta + 2 > \beta + \alpha.$$

Case (ii) follows similarly to (i)(a), we have $\alpha - \epsilon \leq \beta + \epsilon - \delta$,
\begin{eqnarray*}
	\eta(G) & = & p^{\alpha - \epsilon -1}((\beta - \alpha + 2\epsilon - \delta)(p-1) + p+1) \\
	& \geq & 2^{\alpha - \epsilon -1}(\beta - \alpha + 2\epsilon - \delta + 3)\\
	& \geq & \beta - \alpha + 2 \epsilon - \delta + 3(\alpha - \epsilon)\\
	& = & \beta + \alpha + (\alpha - \epsilon - \delta)\\
	& \geq & \beta + \alpha
\end{eqnarray*}
since $\alpha \geq \delta + \epsilon$.$\Box$\\

\section{Metacyclic Groups of Negative Type}

The goal of this section is to compute $\eta$ when $G$ is a metacyclic group of negative type.  We begin by looking at quotients of $G$.  We begin with a preliminary lemma that is useful in understanding the quotients.

Using the notation of Section \ref{secn:meta} and applying Theorem \ref{quot}, we have that if $G = G_2(\alpha, \beta, \epsilon, \delta, - )$ with $\delta \ge 1$ and $N = \langle x^{2^{\alpha - \delta + 1}} \rangle$, then $\eta (G) \ge \eta (G/N)$.  We now show that in fact this is an equality.  We remind the reader that $\alpha - \delta \ge 2$ when $p = 2$.

We now prove the promised equality between $\eta (G)$ and $\eta (G/N)$.  

\begin{theorem}\label{negative quo}
Let $G = G_2 (\alpha, \beta, \epsilon, \delta, - )$ where $\delta \ge 1$.  Then $\eta (G) = \eta (G/N)$ where $N = \langle x^{2^{\alpha-\delta +1}} \rangle$.
\end{theorem}

\noindent
{\bf Proof.}
Note that $N$ does not make sense if $\delta = 0$;  so that it is why we assume $\delta \ge 1$.  Also, if $\delta = 1$, then $N = 1$; so the conclusion is trivial in this case.  Hence, we will assume $\delta \ge 2$.  

We first prove that $\eta (G) = \eta (G/Z)$ where $Z = \langle x^{2^{\alpha-1}} \rangle$.  Recall from Theorem \ref{quot} that to prove $\eta (G) = \eta (G/Z)$, we need to prove that $Z \subseteq G^{\{2\}}$ and if $g \in G \setminus G^{\{2\}}$, then every element of $gZ$ is conjugate to a generator of $\langle g \rangle$.  Observe that $Z \subseteq G^{\{2\}}$.  Since $x^{2^{\alpha - 1}}$ is the only nonidentity element of $Z$, it suffices to prove that if $g \not\in G^{\{2\}}$, then $\langle g \rangle$ and $\langle gx^{2^{\alpha - 1}} \rangle$ are conjugate.  We know from \cite{beuerle} that $G' = \langle x^2 \rangle$.  

We prove the claim by working by induction on $\delta$.  We begin with the case that $\delta = 2$.  We know that $x^y = x^{2^{\alpha - 2} - 1}$.  It follows that 
$$(x^2)^y = (x^y)^2 = (x^{2^{\alpha -2} -1})^2 = x^{2^{\alpha - 1} - 2} = (x^{-2})x^{2^{\alpha - 1}}.$$  
Observe that this yields that $(x^{-2})^y = x^2 x^{2^{\alpha -1}}$.  Using this fact and the observation that $x^{2^{\alpha -1}}$ is central, we then have 
$$(x^2)^{y^2} = (x^{-2} x^{2^{\alpha -1}})^y = x^2 x^{2^{\alpha -1}} x^{2^{\alpha -1}} = x^2.$$  
It follows that $x^2 $ and $y^2$ commute.  Let $A = \langle x^2, y^2 \rangle$, and observe that $G' \le A$, so $A$ is a normal, abelian subgroup of $G$.  

We know that every element of $G$ has the form $y^k x^m$ where $0 \le k \le 2^\beta - 1$ and $0 \le m \le 2^\alpha - 1$ are integers.  Notice that if $4$ divides both $k$ and $m$, then $g \in A^{\{2\}} \subseteq G^{\{2\}}$.  Also, $x^2, y^2 \in G^{\{2\}}$.  

If $g = y^{2l+1}x^m$ for integers $l$ and $m$, then we can appeal to Lemma \ref{conjugacy} to see that $g$ is conjugate to $g x^{2^{\alpha - 1}}$ and so, $\langle g \rangle$ and $\langle gx^{2^{\alpha - 1}} \rangle$ are conjugate, as desired.

Since $x^y = x^{2^{\alpha - 2} - 1}$, we have $x^{y^2} = x^{2^{2\alpha -4} - 2^{\alpha - 1} + 1}$.  Since $\delta \ge 2$, we know that $\alpha \ge 4$ (this is using the fact that $\alpha - \delta \ge 2$), so $2\alpha - 4 \geq  \alpha$.  Hence, we have $x^{y^2} = x^{-2^{\alpha - 1} + 1}$.  In addition, $x^{2^{\alpha - 1}}$ has order $2$, so $x^{-2^{\alpha - 1}} = x^{2^{\alpha - 1}}$.  Thus, we have shown $x^{y^2} = x^{2^{\alpha - 1} + 1}$.

Suppose now that $g = y^{2l} x^{2h +1}$ for integers $l$ and $h$.  From above, we have 
$$g^{y^2} = (y^{2l}x^{2h+1})^{y^2} = y^{2l} (x^{y^2})^{2h+1} = y^{2l}(x^{2^{\alpha -1} + 1})^{2h+1} = y^{2l} x^{2h+1}x^{2^{\alpha - 1}} = g x^{2^{\alpha - 1}}.$$  
We deduce that $\langle g \rangle$ and $\langle gx^{2^{\alpha - 1}} \rangle$ are conjugate, as desired.

We have shown that $x^{y^2} = x x^{2^{\alpha -1}}$.  This implies that $x^{-1} y^{-2} x = x^{2^{\alpha - 1}} y^{-2}$.  Inverting, we obtain $(y^2)^x = y^2 x^{{2^\alpha -1}}$.  Now, suppose that $g = y^{2l}x^{2h}$.  We can assume from above that either $l$ is odd or $h$ is odd.  Assume first that $l$ is odd.  We have 
$$g^{x} = (y^{2l}x^{2h})^x = ((y^2)^x)^l x^{2h} = (y^2 x^{2^{\alpha -1}})^l x^{2h} = y^{2l} x^{2h} x^{2^{\alpha -1}} = g x^{2^\alpha - 1}.$$
We obtain  $\langle g \rangle$ and $\langle gx^{2^{\alpha - 1}} \rangle$ are conjugate, as desired.

We are left with the case that $g = y^{4l} x^{2(2h+1)}$ for integers $h$ and $l$.  We claim that $g \in G^{\{ 2 \}}$.  Notice that there is an integer $k$ so that $\langle g \rangle = \langle y^{4k} x^2 \rangle$ and that $g \in G^{\{ 2\}}$ if and only if $y^{4k} x^2 \in G^{\{2\}}$.  We show that $y^{4k}x^2 \in G^{\{2\}}$.  We have $x^{y^2} = x x^{2^{\alpha -1}}$.  It follows that $xy^2 = y^2 x x^{2^{\alpha - 1}}$ and 
$$(y^{2k}x)^2 = y^{2k}xy^{2k}x = y^{2k} y^{2k}xx^{2^{\alpha - 1}k}x = y^{4k} x^2 x^{2^{\alpha - 1}k}.$$  
When $k$ is even, we see that $(y^{2k}x)^2 = y^{4k} x^2$.  Now assume that $k$ is odd.  We have 
\begin{eqnarray*}
(y^{2k}xx^{2^{\alpha - 2}})^2 &=& y^{2k}x x^{2^{\alpha - 2}}y^{2k}xx^{2^{\alpha - 2}} = y^{2k} y^{2k}xx^{2^{\alpha - 1}k}x x^{2^{\alpha - 2}2}\\
&=& y^{4k} x^2 x^{2^{\alpha - 1}(k+1)} = y^{4k}x^2.
\end{eqnarray*}  
Note that we are using the fact that $x^{2^{\alpha - 2}}$ commutes with both $x$ and $y^2$ here.  Thus, this yields $g \in G^{\{2\}}$.   We conclude for all elements $g \in G \setminus G^{\{2\}}$ that $g$ and $gx^{2^{\alpha -1}}$ are conjugate and we have proved that $\eta (G) = \eta (G/Z)$ when $\delta = 2$.

We now assume that $\delta > 2$.  Let $M = \langle x^2, y \rangle$.  Since $x^y = x^{2^{\alpha - \delta} - 1}$, we see that $(x^2)^y = (x^{2^{\alpha - \delta}-1})^2 = (x^2)^{2^{(\alpha -1) - (\delta - 1)}-1}$.  Also, $y^{2^\beta} = x^{2^{\alpha - \epsilon}} = (x^2)^{2^{(\alpha -1) - \epsilon}}$.  Observe that $(x^2)^{2^{(\alpha - 1) -1}} = x^{2^{\alpha -1}}$.  We conclude that $M = G_2( \alpha - 1, \beta, \epsilon, \delta - 1, -)$.  Let $g \in G \setminus G^{\{2\}}$.  If $g \in M$, then $g \in M \setminus M^{\{2\}}$.   By induction, we have that $g$ is conjugate to $g (x^2)^{2^{(\alpha - 1) - 1}}$, and so, $g$ and $gx^{2^{\alpha - 1}}$ are conjugate.  Thus, we may assume that $g \not\in M$.  This implies that $g = y^l x^{2m + 1}$ for integers $l$ and $m$.  We know that $y$ induces an automorphism of $\langle x \rangle$ of order $2^\delta$.  It follows that $y^{2^{\delta -1}}$ induces an automorphism of $\langle x \rangle$ of order $2$.  Since $\delta \ge 3$, we know that this automorphism is a square.  It is not difficult to see that $x \mapsto x x^{2^{\alpha - 1}}$ is the unique automorphism of $\langle x \rangle$ that has order $2$ and is a square.  Hence, we see that $x^{y^{2^{\delta -1}}} = x x^{2^{\alpha - 1}}$.  We conclude that $g^{y^{2^{\delta -1}}} = (y^l x^{2m+1})^{y^{2^{\delta -1}}} = y^l (x^{y^{2^{\delta -1}}})^{2m+1} = y^l (x x^{2^{\alpha - 1}})^{2m+1} = y^l x^{2m+1} x^{2^{\alpha -1}} = g x^{2^{\alpha - 1}}$.  This completes the proof of the claim that $\eta (G) = \eta (G/Z)$.

We now work to prove $\eta (G) = \eta (G/N)$.  We work by induction on $\delta$.  If $\delta = 2$, then $N = Z$, and the above claim yields the result.  We assume that $\delta \ge 3$.  We have that $\eta (G) = \eta (G/Z)$.  By induction, $\eta (G/Z) = \eta ((G/Z)/(N/Z))$, and the First Isomorphism Theorem implies that $G/N \cong (G/Z)/(N/Z)$, so $\eta (G/N) = \eta ((G/Z)/(N/Z))$, and we have the desired equality.
$\Box$\\

In light of Theorem \ref{negative quo} and Lemma \ref{quo1}, we see that if we can compute $\eta$ for $G_2 (\alpha,\beta,\epsilon,\delta, -)$ when $\delta = 0, 1$, then we can compute $\eta$ for all metacyclic $2$-groups of negative type.  There are a number of cases to consider when $\delta = 0$ or $1$, and then using these cases, we will compute $\eta$ when $\delta \ge 2$.  Recall that the dihedral $2$-groups are the groups of the form $G_2 (\alpha,1,0,0,-)$, the generalized quaternion $2$-groups are of the form $G_2 (\alpha,1,1,0,-)$, and the semi-dihedral groups are of the form $G_2 (\alpha,1,0,1,-)$.  Also, it is known that $G_2 (\alpha, \beta, 1, 0, - )$ and $G_2 (\alpha, \beta, 1, 1, - )$ are isomorphic for all $\alpha \geq 3$ and $\beta \geq 2$.  Since $\delta \le \beta$, it follows that dihedral, generalized quaternion, and semi-dihedral are the only groups of negative type where $\beta = 1$.  

Thus, we need to analyze the negative metacyclic $2$-groups of type 
$$G_2 (\alpha, \beta, \epsilon, \delta, -)$$ 
with $\beta \geq 2$.  We recall a few facts about the classification of such groups. In particular, for negative type $\epsilon$ is either $0$ or $1$ only. Also the parameters satisfy: $\alpha \geq \delta + 2$ and $\beta \geq \delta$ when $\epsilon = 0$ and $\beta \geq \delta + 1$ when $\epsilon = 1$.

When $\delta = 0$ or $1$, there is a particular abelian normal subgroup $M$ of $G$. For this subgroup $M$, we determine which maximal cyclic subgroups of $M$ are maximal in $G$ and how many maximal cyclic subgroups of $G$ lie outside of $M$. This yields the following result.  Recall that $\eta^* (M)$ is the number of $G$-orbits on the $M$-conjugacy classes of maximal cyclic subgroups of $M$.

\begin{proposition} \label{one prop}
Suppose $G$ is $G_2 (\alpha, \beta, \epsilon, \delta, -)$ where $\delta = 0$ or $1$ and $\beta \geq 2$. Let $M = \langle x, y^2 \rangle$. Then $M$ is a normal abelian subgroup of $G$ and the following holds:
\begin{enumerate}
\item[(i)] If $\delta = 0$, then $\eta(G) = \eta^*(M) + 1$ and every maximal cyclic subgroup of $M$ is maximal cyclic in $G$ except $\langle y^2 \rangle$.
\item[(ii)] If $\delta = 1$, then $\eta(G) = \eta^*(M)$ and every maximal cyclic subgroup of $M$ is maximal cyclic in $G$ except $\langle y^2 \rangle$ and $\langle y^2 x^{2^{\alpha -1}} \rangle$.
\end{enumerate}	
\end{proposition}

\noindent
{\bf Proof.} As $M$ is a subgroup of index $2$ in $G$ it follows that $M$ is normal in $G$.  Let $Y = \langle y^2 \rangle$.  Observe that $y^2$ centralizes $\langle x \rangle$ and is obviously central in $\langle y \rangle$; so $Y =  \langle y^2 \rangle$ is central in $G$.  Now, $M$ is central-by-cyclic, so $M$ is abelian.

We now prove that there are exactly two conjugacy classes of maximal cyclic subgroups of $G$ outside of $M$.  Since $\langle x \rangle$ is normal in $G$ and $G = \langle x \rangle \langle y \rangle  = \langle x \rangle \langle xy \rangle$, we see that $C_G (\langle y \rangle) = C_{\langle x \rangle} (\langle y \rangle) \langle y \rangle = \langle x^{2^{\alpha-1}} \rangle \langle y \rangle$ and $C_G (\langle xy \rangle) =  \langle x^{2^{\alpha-1}} \rangle \langle xy \rangle$.  It follows that both $\langle y \rangle$ and $\langle xy \rangle$ lie in conjugacy classes of size $|\langle x \rangle:\langle x^{2^{\alpha-1}}\rangle | = 2^{\alpha-1}$.  It is not difficult to see now that every cyclic subgroup of $G$ outside of $M$ is conjugate to either $\langle y \rangle$ or $\langle xy \rangle$.

$(i)$ For $\delta = 0$ we show that every maximal cyclic subgroup of $M$ is a maximal cyclic subgroup of $G$ except 
$\langle y^2 \rangle$ which lies in exactly $2$ different conjugacy classes of maximal cyclic subgroups of $G$, namely
$\langle y \rangle$ and $\langle xy \rangle$.

Observe that $yY$ acts on $M/Y$ inverting every element.  Thus, $M/Y$
is a cyclic subgroup of index $2$ in $G/Y$. We have $(yY)^2 = Y$, so
$G/Y$ is a dihedral group. It follows that if $g \in G \setminus M$, then
$(gY)^2 = Y$ and so, $g^2 \in Y$. Hence, $Y$ is the only maximal cyclic
subgroup of $M$ that is not maximal cyclic in $G$. Notice that $Y \le \langle y \rangle$. Also, we know that $\langle yY \rangle$ and $\langle xyY \rangle$ are in
different conjugacy classes of subgroups of $G/Y$, so $\langle y \rangle$ and $\langle xy \rangle$ are in different conjugacy classes of $G$. Since $x^y = x^{-1}$, so
$ xy = yx^{-1}$. It follows that $(yx)^2 = yxyx= y(yx^{-1})x = y^2$.

$(ii)$ For $\delta =1$ we show that the only maximal cyclic subgroups of $M$ that are not maximal in $G$ are $\langle y^2 \rangle$
and $\langle y^2 x^{2^{\alpha - 1}} \rangle$. Again there are exactly $2$ different conjugacy classes of maximal 
cyclic subgroups outside of $M$ given by $\langle y \rangle$ and $\langle xy \rangle$. Note that $\langle y \rangle$ contains 
$\langle y^2 \rangle$ and  $\langle xy \rangle$ contains $\langle y^2 x^{2^{\alpha -1}} \rangle$.

Note that
$M/Y$ is cyclic in $G/Y$ of order $2^\alpha$. Also, $(yY)^2 = Y$ and 
$(xY)^{yY} = x^{2^{\alpha-1} - 1} Y = (xY)^{2^{\alpha-1} - 1}$. It follows that $G/Y$ is
isomorphic to a semi-dihedral group. Let $Z = \langle x^{2^{\alpha-1}}, Y \rangle$, and
observe that $Z/Y = Z(G/Y)$. Notice that if $g \in G \setminus M$, then
$(gY)^2 \in Z/Y$. This implies that $g^2 \in Z$. Observe that $\langle y^2 \rangle$ and
$\langle y^2 x^{2^{\alpha-1}} \rangle$ are central (and hence normal) in $G$. It follows
that the square of any conjugate of $y$ will be $y^2$. Since $\delta = 1$, we
have $x^y = x^{2^{\alpha-1} - 1}$, so $xy = yx^{2^{\alpha-1} - 1}$. We have 
$(yx)^2 = yxyx = y(yx^{2^{\alpha-1} - 1})x = y^2 x^{2^{\alpha-1}}$. This implies that the square
of any conjugate of $xy$ will be $y^2 x^{2^{\alpha-1}}$. Hence, any other
subgroup of $M$ that is maximal cyclic in $M$ will be maximal cyclic in $G$.
$\Box$\\[1ex]

We now work to compute $\eta$ for the groups with negative type and $\delta$ equal to $0$ or $1$.  We will first handle the case when $\epsilon = 0$ and $\beta = 2$.  For the following lemma recall that $\alpha \geq \delta +2$ when $p = 2$, so when $\delta =1$ we must have $\alpha \geq 3$.

\begin{lemma} \label{epsilon=0}
Suppose $G$ is $G_2 (\alpha, 2, 0, \delta, - )$. Then 
\begin{enumerate}
	\item[(i)] $\eta(G) = \alpha + 3$ if $\delta = 0$ and
	\item[(ii)] $\eta(G) = \alpha + 2$ if $\delta = 1$.
\end{enumerate}
\end{lemma}

\noindent
{\bf Proof.} 
Following Proposition \ref{one prop}, we take $M = \langle x, y^2 \rangle$; so $M$ is abelian.  We have $M \cong C_{2^{\alpha}} \times C_2$ and $\eta(M) = \alpha +2$ by Lemma \ref{two}.  We claim that all subgroups of $M$ are normal in $G$.  To see this, note that if $K$ is a subgroup of $M$ then (1) $K$ is a subgroup of $\langle x \rangle$, (2) $K = \langle x^a, y^2 \rangle$ for some integer $1 \le a \le 2^\alpha - 1$ or (3) $K = \langle x^a y^2 \rangle$ for some integer $1 \le a \le 2^\alpha -1$.  When $\delta = 0$, we know that $x^y = x^{-1}$, so $(x^a)^y = (x^a)^{-1}$ and $(x^ay^2)^y = (x^ay^2)^{-1}$ for every integer $a$.  When $\delta = 1$, we have $(x^a)^y = x^{a(2^{\alpha - 1} - 1)}$.  The observation is that $\langle x^a \rangle = \langle x^{a(2^{\alpha - 1} - 1)} \rangle$, $\langle x^a, y^2 \rangle = \langle x^{a(2^{\alpha - 1} - 1)}, y^2 \rangle$, and $\langle x^a y^2 \rangle = \langle x^{a(2^{\alpha - 1} - 1)} y^2 \rangle$. This proves the claim.  Therefore $\eta^*(M) = \eta(M)$ and the result follows from Proposition \ref{one prop}.
$\Box$.\\

We continue with the case where $\epsilon = 0$.  We now consider the case that $\beta \ge 3$.  Recall that $g_p (a,b) = p^{(l-1)}((k-l) (p-1) + p +  1)$ where $p$ a prime, $a$ and $b$ are positive integers, and we take $k = {\rm max} (a,b)$ and $l = {\rm min} (a,b)$.  Recall also that $g_p (a,b) = \eta (C_{p^a} \times C_{p^b})$. The following can be viewed as an improvement on Proposition \ref{centre}(ii).

\begin{theorem} \label{large}
Suppose $G$ is $G_2 (\alpha, \beta, 0, \delta, -)$ with $\beta \ge 3$. As previously let $M = \langle x, y^2 \rangle$.
Then the following hold:
\begin{enumerate}
	\item If $\delta = 1$, then $\eta (G) = \eta(M)/2 + 2 = g_2 (\alpha,\beta - 1)/2 + 2$.
	\item If $\delta = 0$, then $\eta (G) =\eta(M)/2 + 3 =  g_2 (\alpha,\beta - 1)/2 + 3$.
\end{enumerate}
\end{theorem}

\noindent
{\bf  Proof.} 
Note that we are assuming $\delta$ is $0$ or $1$. As in Proposition \ref{one prop}, we let $M = \langle x, y^2 \rangle$; so it follows that $M$ is abelian.  In particular, since we are assuming that $\epsilon = 0$, we have $M \cong \langle x \rangle \times \langle y^2 \rangle = C_{2^\alpha} \times C_{2^{\beta - 1}}$.  Using Lemma \ref{two}, we obtain $\eta (M) = g_2 (\alpha,\beta -1)$.  Let $k$ be the maximum of $\alpha$ and $\beta - 1$ and let $l$ be the minimum of $\alpha$ and $\beta - 1$; so that $\eta (M) = g_2 (\alpha,\beta -1) = 2^{l-1} (k-l + 3)$.  We now work to prove that $\eta^* (M) = g_2(\alpha,\beta -1)/2 +2$.   Once this is done, then we will have the conclusion via Proposition \ref{one prop}.

It is not difficult to see that $\langle x \rangle$, $\langle y^2 \rangle$, and $\langle y^2 x^{2^{\alpha - 1}} \rangle$ are maximal cyclic subgroups of $M$ that are normal in $G$.  We claim that $\langle y^{2(2^{\beta -2})} x \rangle$ is a maximal cyclic subgroup of $M$ that is normal in $G$.  It is easy to see that it is maximal cyclic.  When $\delta = 0$, we see that $(\langle y^{2(2^{\beta -2})} x \rangle)^y = \langle y^{2(2^{\beta -2})} x^{-1} \rangle = \langle (y^{2(2^{\beta -2})} x)^{-1} \rangle$, and when $\delta = 1$, we have $(\langle y^{2(2^{\beta -2})} x \rangle)^y = \langle y^{2(2^{\beta -2})} x^{2^{\alpha}-1} \rangle = \langle (y^{2(2^{\beta -2})} x)^{2^{\alpha}-1} \rangle$.  This proves that it is normal in $G$.

We will prove that all the other maximal cyclic subgroups of $M$ will be in conjugacy classes of size $2$ in $G$.  Thus, $\eta^*(M) = (\eta(M) - 4)/2 + 4 = \eta (M)/2 -2 +4 =  g_2 (\alpha,\beta - 1)/2 + 2$.

Let $C$ be a maximal cyclic subgroup of $M$.  It is not difficult to see that $C$ will be generated by an element of the form $y^{2l}x$ or one of the form $y^2 x^l$.  When $\delta = 0$, we have that $(y^{2l}x)^y = y^{2l}x^{-1}$ and $(y^2x^l)^y = y^2 x^{-l}$.  For $C$ to be normal, we need this conjugate to be in $C$.  When the generator is $y^{2l}x$, we need $y^{2l}x^{-1} = (y^{2l}x)^k = y^{2lk}x^k$ for some integer $k$.  This implies that $y^{2l-2lk} = x^{k+1}$.  Since $\epsilon = 0$, we have that $y^{2l-2lk} = x^{k+1} = 1$.  We see that we must have $2^\alpha$ dividing $k+1$ and $2^{\beta}$ must divide $2l(1-k)$.  Thus, there is an integer $r$ so that $k+1 = 2^\alpha r$, and thus, $k = 2^\alpha r - 1$.  We obtain that $2^{\beta -1}$ must divide $l(1-(2^\alpha r -1)) = l(2 - 2^\alpha r) = 2l (1 - 2^{\alpha - 1}r)$.  Since we know that $\alpha \ge 2$, this implies that $2^{\beta - 2}$ must divide $l$.  It follows that $\langle x \rangle$ and $\langle y^{2^{\beta-1}} x \rangle$ are the only two maximal cyclic subgroups of $M$ that are normal in $G$ that are generated by an element of the form $y^{2l}x$ when $\delta = 0$. 

When the generator is $y^2 x^l$, we need $y^2x^{-l} = (y^2 x^l)^k = y^{2k}x^{lk}$ for some integer $k$.  This implies that $y^{2-2k} = x^{lk+l} = 1$.  This implies that $2^\beta$ divides $2(1-k)$ and so, $2^{\beta -1}$ divides $1-k$.  Hence, there is an integer $r$ so that $1-k = r2^{\beta -1}$, and hence, $k = 1 -r2^{\beta -1}$.  We see that $2^\alpha$ divides $l(1 + k) = l (1 + (1 - r2^{\beta -1})) = l (2-r2^{\beta-1}) = 2l(1-r2^{\beta - 2})$.  Since $\beta \ge 3$, we deduce that $2^{\alpha - 1}$ must divide $l$.  It follows that $\langle y^2 \rangle$ and $\langle y^2 x^{2^{\alpha - 1}} \rangle$ are the only maximal cyclic subgroups of $M$ that are normal in $G$ that are generated by an element of the form $y^2 x^l$ when $\delta = 0$.  This proves the result when $\delta = 0$.

Now we suppose that $\delta = 1$.  Recall that $\alpha \ge \delta + 2$, so $\alpha \ge 3$.  We have that $(y^{2l}x)^y = y^{2l}x^{2^{\alpha-1}-1}$ and $(y^2x^l)^y = y^2 x^{l(2^{\alpha-1}-1)}$.  For $C$ to be normal, we need this conjugate to be in $C$.  Suppose the generator is $y^{2l}x$. We need $y^{2l}x^{2^{\alpha-1}-1} = (y^{2l}x)^k = y^{2lk}x^k$ for some integer $k$.  This implies that $y^{2l-2lk} = x^{k-2^{\alpha -1} +1} = 1$.  We deduce that $2^\alpha$ must divide $k - 2^{\alpha -1} + 1$, and so, there is an integer $r$ so that $k - 2^{\alpha -1} + 1 = 2^\alpha r$.  We obtain $k = 2^\alpha r + 2^{\alpha - 1} - 1$.  We have that $2^\beta$ divides $2l(1-k) = 2l(1 - 2^\alpha r - 2^{\alpha - 1} + 1)$.  It follows that $2^{\beta - 2}$ divides $l(1 - 2^{\alpha -1}r -2^{\alpha - 2})$.  Since $\alpha \ge 3$, we see that $2^{\beta - 2}$ divides $l$.  We conclude that $\langle x \rangle$ and $\langle y^{2^{\beta-1}} x \rangle$ are the only two maximal cyclic subgroups of $M$ that are normal in $G$ that are generated by an element of the form $y^{2l}x$ when $\delta = 1$.

When the generator is $y^2 x^l$, we need $y^2x^{l(2^{\alpha - 1}-1)} = (y^2 x^l)^k = y^{2k}x^{lk}$ for some integer $k$.  We see that $y^{2-2k} = x^{lk - l(2^{\alpha -1} - 1)} = 1$.   It follows that $2^\beta$ divides $2 (1-k)$, and so, $2^{\beta -1}$ divides $1-k$.  There is an integer $r$ so that $1-k = 2^{\beta -1}r$ which yields $k = 1 - 2^{\beta - 1}r$.  We now determine that $2^\alpha$ divides $l(k - 2^{\alpha -1} + 1) = l (1 - 2^{\beta -1}r -2^{\alpha -1} + 1) = 2l (1 - 2^{\beta - 2}r -2^{\alpha - 2})$.  Since $\alpha \ge 3$ and $\beta \ge 3$, we have that $2^{\alpha -1}$ divides $l$.  We conclude that $\langle y^2 \rangle$ and $\langle y^2 x^{2^{\alpha - 1}} \rangle$ are the only maximal cyclic subgroups of $M$ that are normal in $G$ that are generated by an element of the form $y^2 x^l$ when $\delta = 1$.  This proves the result when $\delta = 1$.
$\Box$\\

In this next corollary, recall that $\delta \le \beta$, so when $\beta = 2$, we must have $\delta = 2$.  We are able to use Theorem \ref{large} to compute $\eta$ for groups of negative type where $\delta \ge 2$.

\begin{corollary} \label{delta ge 2}
Suppose $G$ is $G_2 (\alpha,\beta,\epsilon,\delta,-)$ with $\delta \ge 2$, then 
\begin{enumerate}
	\item $\eta (G) = \alpha - \delta + 3 = \alpha + 1$ if $\beta = 2$.
	\item $\eta (G) = g_2 (\alpha - \delta +1,\beta-1)/2 + 2$ if $\beta \ge 3$.
\end{enumerate}
\end{corollary}

\noindent
{\bf Proof.}
By Theorem \ref{negative quo}, we have that $\eta (G) = \eta (G/N)$ where $N = \langle x^{2^{\alpha - \delta + 1}} \rangle$.  Applying Lemma \ref{quo1}, we see that $G/N \cong G_2 (\alpha - \delta +1, \beta, 0,1, -)$.  Using Lemma \ref{epsilon=0}, we see that $\eta (G/N) = \alpha - \delta + 1 + 2 = \alpha + 3 - \delta$ when $\beta = 2$.  Since $2 \leq \delta \leq \beta =2$, we see that $\delta = 2$, and so, $\eta (G) = \alpha + 1$.  When $\beta \ge 3$, we apply Theorem \ref{large} to see that $\eta (G) = \eta (G/N) = g_2 (\alpha - \delta +1, \beta -1)/2 + 2$.
$\Box$\\

We now compute $\eta$ for groups of negative type with $\delta = 0$ and $\epsilon = 1$.  We first handle the case where $\beta = 2$.  

\begin{lemma} \label{beta=2}
Suppose $G$ is $G_2 (\alpha, 2, 1, 0, -)$ then $\eta(G) = \alpha + 2$.
\end{lemma}

\noindent
{\bf Proof.}  
Define $M = \langle x, y^2 \rangle$.  By Proposition \ref{one prop}, we know that $M$ is a normal abelian subgroup of $G$.  First note that $(x^{2^{\alpha - 2}}y^2)^2 = x^{2^{\alpha-1}}y^4 = x^{2^{\alpha -1}}x^{2^{\alpha -1}} = x^{2^{\alpha}} =1.$ Thus, $M = \langle x \rangle \times \langle x^{2^{\alpha - 2}} y^2 \rangle \cong C_{2^\alpha} \times C_2$ and $\eta(M) = \alpha + 2$. Consideration of the maximal cyclic subgroups of $M$ shows that all are normal except $\langle (1, x^{2^{\alpha -2}}y^2) \rangle$ and $\langle (x^{2^{\alpha-1}}, x^{2^{\alpha - 2}}y^2) \rangle$ which are conjugate in $G$ via $y$.  To see that these two subgroups are conjugate, observe that $M$ has three subgroups of order $2$ and that $\langle x^{2^{\alpha - 1}} \rangle = \langle y^{2^\beta} \rangle$ is central in $G$ and that $Z(G)$ is cyclic.  Either $y$ normalizes both of the other two subgroups of order $2$ or it permutes them.  However, if $y$ were to normalize them, they would be normal in $G$ and since they have order $2$, that would imply that they would be central in $G$.  This however would contradict the fact that the center of $G$ is cyclic. Thus $\eta^*(M) = \alpha + 1$. The result follows from Proposition \ref{one prop}. 
$\Box$\\

We continue with the groups of negative type where $\delta = 0$ and $\epsilon = 1$.  We next consider $\beta \ge 3$ and $\alpha = 2$. 

\begin{lemma} \label{betage3}
Suppose $G$ is $G_2 ( 2, \beta, 1, 0 , - )$ with $\beta \geq 3$. Then $\eta(G) = \beta + 2$.
\end{lemma}

\noindent
{\bf Proof.} 
Define $M = \langle x, y^2 \rangle$.  By Proposition \ref{one prop}, we know that $M$ is a normal abelian subgroup of $G$.  Note that $(xy^{2^{\beta -1}})^2 = x^2y^{2^{\beta}} = x^2x^2 = x^4 = 1$. So $M = \langle xy^{2^{\beta -1}} \rangle \times \langle y^2 \rangle \cong C_2 \times C_{2^{\beta}}$ and $\eta(M) = \beta + 2$. Consideration of the maximal cyclic subgroups of $M$ shows that all are normal except for $\langle (xy^{2^{\beta -1}}, 1) \rangle$ and $\langle (xy^{2^{\beta -1}}, y^{2^{\beta}}) \rangle$ which are conjugate in $G$ via $y$.  The proof that these two subgroups are conjugate is similar to the proof of Lemma \ref{beta=2}.  In particular, $Z(G)$ is cyclic, $M$ has three subgroups of order $2$, and if $y$ normalized these two subgroups, then it would centralize them and contradict the fact that $Z (G)$ is cyclic. Thus $\eta^*(M) = \beta + 1$.  The result follows from Proposition \ref{one prop}.
$\Box$\\

We conclude by computing $\eta$ when $\delta = 0$, $\epsilon = 1$, $\alpha \ge 3$, and $\beta \ge 3$. Note this also covers the cases $\delta =1$,
$\epsilon =1$ and $\alpha, \beta \ge 3$.

\begin{theorem} \label{largeone}
Suppose $G$ is $G_2 (\alpha, \beta, 1, 0, -)$ with $\alpha \ge 3$ and $\beta \ge 3$.  Let $M = \langle x, y^2 \rangle$.
\begin{enumerate}
	\item If $\alpha \ge \beta$, then $\eta (G) = \eta(M)/2 + 3 = g_2(\alpha,\beta-1)/2 + 3$.
	\item If $\alpha < \beta$, then $\eta (G) = \eta(M)/2+ 3 =  g_2 (\alpha-1,\beta)/2 + 3$.
\end{enumerate}
\end{theorem}

\noindent
{\bf Proof.}
As in Proposition \ref{one prop}, we let $M = \langle x, y^2 \rangle$; so it follows that $M$ is abelian.  We know that $|M| = 2^{\alpha + \beta - 1}$, that $x$ has order $2^\alpha$ and $y^2$ has order $2^\beta$.  Suppose $\alpha \ge \beta$, then $M \cong C_{2^\alpha} \times C_{2^{\beta -1}}$, and so $\eta (M) = g_2 (\alpha,\beta-1)$.  Let $w = y^2 x^{2^{\alpha - \beta}}$.  Observe that $w^{2^{\beta-2}} = (y^2x^{2^{\alpha-\beta}})^{2^{\beta-2}} = y^{2^{\beta -1}} x^{2^{\alpha - 2}} \not\in \langle x \rangle$ and $w^{2^{\beta-1}} = (y^2x^{2^{\alpha-\beta}})^{2^{\beta-1}} = y^{2^{\beta}} x^{2^{\alpha - 1}} = x^{2^{\alpha -1}} x^{2^{\alpha -1}} = 1$.  It follows that $M = \langle x \rangle \times \langle w \rangle$.  

If $\beta \ge \alpha + 1$, then $M \cong C_{2^{\alpha-1}} \times C_{2^{\beta}}$, and so $\eta (M) = g_2 (\alpha-1,\beta)$.  Let $u = y^{2^{\beta-\alpha+1}}x$.  We compute $u^{2^{\alpha -2}} = (y^{2^{\beta-\alpha+1}}x)^{2^{\alpha -2}}= y^{2^{\beta-1}}x^{2^{\alpha - 2}} \not\in \langle y \rangle$ and $u^{2^{\alpha -1}} = (y^{2^{\beta-\alpha+1}}x)^{2^{\alpha -1}}= y^{2^{\beta}}x^{2^{\alpha - 1}} = x^{2^{\alpha -1}} x^{2^{\alpha -1}} = 1$.  We deduce that $M = \langle u \rangle \times \langle y \rangle$.

In both cases, we will show that $\eta^* (M) = \eta (M)/2 + 2$, and we obtain the conclusion by applying Proposition \ref{one prop}.  Notice that a maximal cyclic subgroup of $M$ will be generated either by an element of the form $y^{2l}x$ for some integer $l$ or by an element of the form $y^2 x^l$ for some integer $l$.  Observe that $\langle x \rangle$ and $\langle y^2 \rangle$ are maximal cyclic subgroups of $M$ that are normal in $G$.  

We next show that $\langle y^{2^{\beta -1}} x\rangle$ and $\langle y^2 x^{2^{\alpha -2}} \rangle$ are normal subgroups in $G$. Since $M$ is abelian and has index $2$ in $G$, it suffices to show that $y$ normalizes these subgroups. We compute $(y^{2^{\beta -1}} x)^y = y^{2^{\beta -1}} x^{-1} = (y^{2^{\beta -1}} x)^{-1}$.  Since $y$ conjugates the generator of $\langle y^{2^{\beta-1}}x \rangle $ to its inverse, this implies that $\langle y^{2^{\beta -1}} x\rangle$ is normal in $G$.   

We now turn to $\langle y^2 x^{2^{\alpha -2}} \rangle$.  We begin with the observation that $(y^2 x^{2^{\alpha - 2}})^4 = y^8$.  Since $\beta \ge 3$, we see that $x^{2^{\alpha -1}} = y^{2^\beta} \in \langle y^2 x^{2^{\alpha -2}} \rangle$.   Conjugating yields $(y^2 x^{2^{\alpha -2}})^y = y^2 x^{-2^{\alpha - 2}}$.  Note that $x^{-2^{\alpha -2}} = x^{2^{\alpha -2}} x^{2^{\alpha -1}}$.  We have $(y^2 x^{2^{\alpha -2}})^y = y^2 x^{2^{\alpha -2}} x^{2^{\alpha -1}}$.  Since both $y^2 x^{2^{\alpha -2}}$ and $x^{2^{\alpha -1}}$ lie in $\langle y^2 x^{2^{\alpha -2}} \rangle$, we conclude that $(y^2 x^{2^{\alpha -2}})^y$ lies in $\langle y^2 x^{2^{\alpha -2}} \rangle$.
We deduce that $\langle y^2 x^{2^{\alpha -2}} \rangle$ is normal in $G$. 

We prove that the remaining maximal cyclic subgroups of $M$ lie in orbits of size $2$.  We have noted that a maximal cyclic subgroup $C$ of $M$ will have a generator of the form $y^{2l} x$ or of the form $y^2 x^l$ for some integer $l$.  If $C$ has a generator of the form $y^{2l} x$, then for $C$ to be normal we need $(y^{2l}x)^y = y^{2l}x^{-1} \in C$.  This implies that $y^{2l} x^{-1} = (y^{2l}x)^k$ for some integer $k$.  We have $y^{2l-2lk} = x^{k+1} = u \in \langle x \rangle \cap \langle y^2 \rangle = \langle x^{2^{\alpha - 1}} \rangle$.  Hence, $u$ is either $1$ or $x^{2^{\alpha - 1}}$.  If $u = 1$, then $2^{\alpha}$ divides $k+1$ and $2^{\beta +1}$ divides $2l(1-k)$.  We see that there is an integer $r$ so that $k+1 = 2^{\alpha} r$, and hence, $k = 2^\alpha r - 1$.  This implies that $2^{\beta +1}$ divides $2l (1-k) = 2l(1-2^\alpha r +1)  = 4 l(1 -2^{\alpha -1}r)$.  Since $\alpha \ge 2$, this yields $2^{\beta -1}$ divides $l$.  When $u = x^{2^{\alpha -1}}$, we obtain that $k + 1 \equiv 2^{\alpha -1} ~({\rm mod}~2^{\alpha})$.  Hence, there is an integer $r$ so that $k + 1= 2^{\alpha-1} + r 2^{\alpha}$, and so, $k = 2^{\alpha - 1} + r2^\alpha - 1$.  We see that $2l(1-k) \equiv 2^\beta ~({\rm mod}~2^{\beta+1})$.  This implies that $2^{\beta+1}$ divides $2l(1-k) - 2^{\beta} = 2l (1 -2^{\alpha -1} - r 2^\alpha + 1) -2^\beta = 4l(1-2^{\alpha -2} -r 2^{\alpha-1}) -2^\beta$.  We deduce that $2^{\beta-2}$ divides $l$.  We conclude that $\langle x \rangle$ and $\langle y^{2^{\beta -1}} x \rangle$ are the only maximal cyclic subgroups of $M$ having the form $\langle y^{2l}x \rangle$ that are normal in $G$.

We now suppose that $C$ has a generator of the form $y^2 x^l$.  We need $(y^2 x^l)^y = y^2 x^{-l} \in C$.  Hence, we have that $y^2 x^{-l} = (y^2 x^l)^k = y^{2k} x^{lk}$ for some integer $k$.  We have $y^{2-2k} = x^{lk+l} = u$.  As in the previous paragraph, we see that $u$ is either $1$ or $x^{2^{\alpha -1}}$.  If $u = 1$, then we have that $2^{\beta +1}$ divides $2(1-k)$, and so, there is an integer $r$ so that $1-k = 2^\beta r$.  We determine that $2^\alpha$ divides $l(k+1) = l(1 - 2^\beta r +1) = 2l(1-2^{\beta-1}r)$.  It follows that $2^{\alpha - 1}$ divides $l$.  Now, suppose that $u = x^{2^{\alpha -1}}$.  We must have that $2(1-k) \equiv 2^\beta ~({\rm mod}~2^{\beta +1})$ and $l(k+1) \equiv 2^{\alpha -1} ~({\rm mod}~2^\alpha)$.  Hence, there is an integer $r$ so that $2(1-k) = 2^\beta + 2^{\beta + 1}r$.  This implies that $k = 1 -2^{\beta -1} - 2^{\beta} r$.  We then obtain that $2^\alpha$ divides $l (k+1) - 2^{\alpha - 1} = l (1 -2^{\beta-1} - 2^\beta r + 1) -2^{\alpha -1} = 2 (l(1 - 2^{\beta -2} - 2^{\beta -1}r) - 2^{\alpha - 2})$.  This implies that $2^{\alpha -1}$ divides $l (1 - 2^{\beta -2} - 2^{\beta -1}r) -2^{\alpha -2}$.  Hence, there is an integer $s$ so that $l (1 - 2^{\beta -2} - 2^{\beta -1}r) -2^{\alpha -2} = 2^{\alpha -1}s$.  This leads to $l (1 - 2^{\beta -2} - 2^{\beta -1}r) = 2^{\alpha -1}s + 2^{\alpha -2} = 2^{\alpha -2} (2s + 1)$.  This yields $2^{\alpha -2}$ divides $l$.   Observe that $x^{2^{\alpha - 1}} = y^{2^\beta}$, and so, $\langle y^2 x^{2^{\alpha -1}} \rangle = \langle y^2 \rangle$.  We deduce that $\langle y^2 \rangle$ and $\langle y^2 x^{2^{\alpha - 2}} \rangle$ are the only maximal cyclic subgroups of $M$ having the form $\langle y^2 x^l \rangle$ that are normal in $G$,

We now see that the number of $G$-orbits of maximal cyclic subgroups of $M$ is $(\eta (M) - 4)/2 +4 = \eta (M)/2 - 2 +4 = \eta (M) + 2$, which completes the proof of the result. 
$\Box$\\

We close by proving that when $G$ is metacyclic of minus type that is not dihedral, generalized quaternion, or semi-dihedral, then $\eta (G) \ge \alpha + \beta - 2$ and we determine when equality occurs.  We first handle when $\delta$ equals $0$ or $1$.  In this case, we have $\eta (G) \ge \alpha + \beta$.

\begin{proposition} 
Suppose $G = G_2(\alpha, \beta, \epsilon, \delta, -)$ with  $\delta = 0$ or 1 and $\beta \geq 2$. Then $\eta(G) \geq \alpha + \beta.$
\end{proposition}

\noindent{\bf Proof.}  
(i) Suppose $\epsilon = 0$. 
Denote $l = \min  (\alpha, \beta - 1)$ and $k = \max (\alpha, \beta -1)$.
First, consider $l \geq 3$. Then $\beta \geq 4$ and by Theorem \ref{large} and  Lemma \ref{g2 comp}
$$\eta(G) \geq  g_2(\alpha, \beta -1)/2 + 2 \geq 2k + 2 \geq \alpha + \beta.$$

Next, assume $l = 2$. So $\beta \geq 3$ and by Theorem \ref{large} and Lemma \ref{g2 comp}
$$\eta(G) \ge   g_2(\alpha, \beta -1 )/2 + 2 = k+3 \geq \alpha + \beta.$$

Finally, set $l=1$. As $\alpha \geq 2$, we have $\beta = 2$. The result follows from Lemma \ref{epsilon=0}.

(ii) Now suppose $\epsilon = 1$. Assume  $\alpha \geq \beta$, then $l = \min (\alpha, \beta -1) = \beta -1$
and $k = \max (\alpha, \beta - 1) = \alpha$. If $l \geq 3$, then $\beta \geq 4$ and $\alpha \geq 4$, so we can assume $\delta = 0$. Applying Theorem \ref{largeone} and Lemma \ref{g2 comp} yields
$$\eta(G) =  g_2(\alpha, \beta -1)/2 + 3 \geq 2k + 3 \geq \alpha + \beta.$$
If $l=2$, then $\beta = 3$,
and we again appeal to Theorem \ref{largeone} to obtain 
$$\eta(G) = g_2(\alpha, \beta -1)/2 + 3 = g_2(k, 2)/2 + 3 = k+4 \geq \alpha + \beta.$$
If $l=1$, then $\beta = 2$. If $\alpha = 2$ then $\delta =0$ and if $\alpha \geq 3$ we can assume $\delta =0$.  Thus we apply Lemma \ref{beta=2}. 

Finally, suppose $\epsilon =1$ and $\alpha < \beta$.  
We set $l = \min (\alpha-1, \beta) = \alpha -1$ and $k= \max (\alpha - 1, \beta) = \beta$.
When $l \geq 3$, we apply Theorem \ref{largeone} and Lemma \ref{g2 comp} to get
$$\eta(G) =  g_2(\alpha -1, \beta)/2 + 3 \geq 2k + 3 \geq \alpha + \beta.$$
If $l=2$, then $\alpha = 3$ and $\beta > 3$. Apply  Theorem \ref{largeone} with Lemma \ref{g2 comp} to give
$$\eta(G) =  g_2(\beta, 2)/2 + 3 = \beta +4 \geq \alpha + \beta.$$

If $l=1$, then $\alpha = 2$ and $\delta =0$,  the result follows from Lemma \ref{betage3}. 
$\Box$\\

We now have the case where $\delta \ge 2$.

\begin{proposition} 
Suppose $G = G_2(\alpha, \beta, \epsilon, \delta, -)$ with $\delta \geq 2$.  Then $\eta(G) \geq \alpha + \beta -2$. Equality holds if and only if $\beta = \delta$ and either (i) $\beta=3$ or (ii) $\beta \geq 4$ and $\alpha - \beta =2$. 
\end{proposition}

\noindent{\bf Proof.} 
Set  $l = \min (\alpha - \delta +1, \beta -1)$ and $k = \max (\alpha - \delta +1, \beta-1)$.
We consider various cases according to the value of $l$. 

First, suppose $l \geq 4$. 
Then by Corollary \ref{delta ge 2} and Lemma \ref{g2 comp} 
\begin{eqnarray*}
	\eta(G)      
	& = & g_2(\alpha - \delta +1, \beta-1)/2 + 2 \\
	& = & g_2( k,l)/2 + 2 \geq 2k + l + 2 \\
        &  = & (k+l) + k + 2\\
	& \geq & \alpha - \delta + \beta + \beta -1 + 2 \\
	& \geq & \alpha + \beta + 1\\
\end{eqnarray*} 
since $\delta \leq \beta$. 

Now consider $l = 3$.  We use Corollary \ref{delta ge 2} and Lemma \ref{g2 comp} to find an exact value for $\eta(G)$.
$$\eta(G)     
=   g_2(\alpha - \delta + 1, \beta -1)/2 + 2 
=  g_2( k,3)/2 + 2 
=  2k + 2. $$
If $\alpha -\delta + 1 > \beta-1 = 3$, then $\delta \leq 4$ and
$$\eta(G) = 2(\alpha -  \delta + 1) + 2 = \alpha + (\alpha - \delta + 2) + (-\delta + 2) > \alpha + \beta -2.$$
On the other hand, when $\beta - 1 \geq \alpha - \delta + 1 = 3$, we obtain $\beta \ge 4$ and $\alpha - \delta = 2$, so $\alpha - 2 \leq \beta$ and
$$\eta(G) = 2(\beta - 1) + 2 = 2 \beta \geq \beta + \alpha - 2$$ 
with equality if and only if $\beta = \delta$.

Next suppose $l=2$.  Since $\alpha - \delta + 1 \ge 2+1 = 3$, we must have $\beta =3$.  Applying Corollary \ref{delta ge 2} and Lemma \ref{g2 comp},
\begin{eqnarray*}
\eta(G) & = & g_2(\alpha - \delta +1, \beta - 1)/2 + 2 =  g_2(k, 2)/2 + 2 \\
           & = &  k+3 
             =  \alpha - \delta + 4 \\
           & \geq & \alpha + 1 
             =  \alpha + \beta -2
\end{eqnarray*}
with equality if and only if $\delta = 3 = \beta.$

Lastly consider $l=1$. In this case $\beta = 2$ and the result follows from Corollary \ref{delta ge 2}.
$\Box$\\

\noindent Mariagrazia Bianchi:
Dipartimento di Matematica F. Enriques,
Universit\`a degli Studi di Milano, via Saldini 50,
20133 Milano, Italy.\\
mariagrazia.bianchi@unimi.it\\[1ex]
Rachel Camina: Fitzwilliam College, Cambridge, CB3 0DG, UK.\\
rdc26@cam.ac.uk\\[1ex]
Mark L. Lewis:  Department of Mathematical Sciences, Kent State University, Kent, Ohio, 44242 USA.\\
lewis@math.kent.edu\\[1ex]

\end{document}